\input epsf
\documentclass[12pt]{article}
\usepackage{epsfig}
\usepackage{amssymb}
\usepackage{graphicx}
\newcommand{\be}{\begin{equation}}
\newcommand{\ee}{\end{equation}}
\newcommand{\qed}{\hfill $\square$\vskip .2cm}

\newtheorem{remark}{Remark}[section]

\renewcommand{\S}{{\mathbb S}}

\renewcommand{\limsup}{{\overline{\lim}}}
\renewcommand{\liminf}{{\underline{\lim}}}

\def\<{\langle}
\def\>{\rangle}

\newtheorem{prop}{Proposition}[section]
\newtheorem{defn}{Definition}[section]
\newtheorem{cor}{Corollary}[section]
\newtheorem{lem}{Lemma}[section]
\newtheorem{thm}{Theorem}[section]

\newcommand{\sect}[1]{\section{#1}\setcounter{equation}{0}}

\begin{document}
\title{
Growth of preferential attachment random graphs via continuous-time
branching processes}

\author{Krishna B. Athreya$^1$,\ Arka P. Ghosh$^2$, \ and \  Sunder Sethuraman$^3$}

\thispagestyle{empty}

 \maketitle
 \abstract{
A version of ``preferential attachment'' random
graphs, corresponding to linear ``weights'' with random
 ``edge additions,'' which generalizes some previously considered models, is
 studied. This graph model is embedded in a continuous-time branching
 scheme and, using the branching process apparatus, several
 results
 on the graph model asymptotics are obtained, some extending previous results,
such as
 growth rates for a typical degree and the maximal degree, behavior of the
 vertex where the maximal degree is attained, and a law of
 large numbers for the empirical distribution of degrees which shows certain
 ``scale-free'' or ``power-law'' behaviors.

}

 \vskip .1cm
 \thanks{Research supported in part by NSA-H982300510041, NSF-DMS-0504193
 and NSF-DMS-0608669.  
 \\
 {\sl Abbreviated title}: Preferential attachment random graphs via branching processes  \\
 {\sl AMS (2000) subject classifications}: Primary 05C80; secondary 60J85
 .}
\vskip .1cm

\noindent $^1$ Departments of Mathematics and Statistics, Iowa State University,
Ames, IA \ 50011

\noindent $^2$ Department of Statistics, Iowa State University, Ames, IA \
50011

\noindent $^3$ Department of Mathematics, Iowa State University,
 Ames, IA \ 50011

\noindent {\sl Email:} K.B. Athreya (kba@iastate.edu); A.P. Ghosh
(apghosh@iastate.edu); S. Sethuraman (sethuram@iastate.edu)


\sect{Introduction and results} Preferential attachment processes
have a long history dating back at least to Yule \cite{Yule} and Simon
\cite{Simon} (cf. \cite{Mitzenmacher} for an interesting survey).
Recently, Barabasi and Albert
\cite{BA} proposed a random graph version of these processes as a model for several
real-world networks, such as the internet and various communication
structures, on which there has been much renewed study (see
\cite{Albert_Barabasi}, \cite{CL},
\cite{Durrett_book}, \cite{Newman} and references therein).  To
summarize, the basic idea is that,
starting from a small number of nodes, or vertices, one builds an evolving
graph by ``preferential attachment,'' that is by attaching
new vertices to existing nodes with probabilities proportional to
their ``weight.'' When the weights are increasing functions of the ``connectivity,''
already well connected vertices tend to become even more connected
as time progresses, and so, these graphs can be viewed as types of
``reinforcement'' schemes (cf. \cite{Pemantle}).  A key point, which
makes these graph models ``practical,'' is
that, when the weights are linear, the long term degree proportions
are often in the form of a
``power-law'' distribution whose exponent, by varying parameters, can be matched to
empirical
network data.

The purpose of this note is to understand a
general form of the linear weights model with certain random
``edge additions'' (described below in subsection 1.1)
in terms of an embedding in continuous-time branching processes which allows for
extensions of law of large numbers and maximal degree growth
asymptotics, first approached by difference equations and martingale
methods, in \cite{Bollobas}, \cite{CF}, \cite{Mori1}, \cite{Mori}.

We remark some connections to branching and continuous-time Markov
processes have also been studied in two recent papers.  In
\cite{RTV}, certain laws of large numbers for the degree
distributions of the whole tree, and as seen from a randomly
selected vertex are proved for a class of ``non-explosive'' weights
including linear weights.  In
\cite{OS}, asymptotic degree distributions under super-linear
weights are considered. In this context, the embedding given here
is of a different character with respect to Markov branching systems with immigration, and the contributions made are also
different, concentrating on detailed investigations of a generalized linear weights degree landscape.

\subsection{ Model} Start with two vertices $v_1$, $v_2$ and one edge
joining them--denote this graph as $G_0$. To obtain $G_1$, create a
new vertex $v_3$, and join it a random number $X_1$ times to one of
$v_1$ and $v_2$ of $G_0$ with equal probability. For any finite
graph $G_n= \{v_1,v_2, \ldots, v_{n+2}\}$, let the degree of each
vertex be defined as the number of edges emanating from that vertex,
and the degree of the $j$th vertex, $v_j\in G_n$ be denoted by
$d_j(n)$ for $j=1,\ldots, n+2$ and $n\geq 0$ (note that in our notation, $G_n$ has
$n+2$ vertices at step $n\ge 0$). After $n+2$ vertices are created, to
obtain $G_{n+1}$ from $G_n$, create an $(n+3)$rd vertex $v_{n+3}$,
and connect it a random number $X_{n+1}$ times to one of the $n+2$
existing vertices $v_1,\ldots, v_{n+2}$ with probability
\begin{equation}
\label{weight_eq}\frac{ d_i(n) +\beta}{\sum_{j=1}^{n+2}(
   d_j(n) +\beta)}\end{equation}
of being joined to vertex $v_i$ for $1\leq i\leq n+2$ where
$\beta\geq 0$ is a parameter.
We will also assume throughout that $\{X_i\}_{i\geq 1}$ are
independently and identically distributed  positive integer valued
random variables with distribution $\{p_j\}_{j\geq 1}$ with finite
mean.
The ``weight'' then of the $i$th vertex at the $n$th step is
proportional to
$d_i(n) +\beta$, and linear in the degree.

We remark this basic model creates a growing ``tree'' with
undirected edges, and includes the ``one-edge'' case of the original
Barabasi-Albert process, made precise in \cite{Bollobas}, by setting
$X_i\equiv 1$ and $\beta=0$, as well as the ``$\beta\geq 0$'' scheme
considered in \cite{Mori1} and \cite{Mori}, by taking $X_i\equiv 1$.
Also, the ``$\beta \geq 1$'' linear case considered in \cite{RTV} is
recovered by taking $X_i\equiv 1$.

The aspect of adding a
random number of edges $\{X_i\}_{i\geq 1}$ at each step to vertices
chosen preferentially seems to be a new twist on the standard model
which can be interpreted in various ways.
 The results, as will be seen, involve the mean
 number $\sum jp_j$ of added edges, indicating a sort of ``averaging'' effect
 in the asymptotics.

We also note, in the case $\beta=0$, a more general graph process,
allowing cycles and self-loops, can be formed in terms of the ``tree''
model above (cf. \cite{Bollobas} and
Ch. 4 \cite{Durrett_book})
where several
sets of edges are added to possibly different existing vertices at
each step preferentially.
  Namely, let $\{L_i\}_{i\geq 1}$ be independent and identically
distributed positive integer valued random variables with distribution
$\{q_j\}_{j\geq 1}$ with finite mean, and let $\bar{L}_i =
\sum_{k=1}^iL_k$ for $i\geq 1$.
As before, initially, we start
with two vertices, $v^{(L)}_1$ and $v^{(L)}_2$ and one edge between them.
Run the ``tree'' model now to obtain vertices $\{w_i\}_{i\geq 3}$ and
identify
sets $$
\{w_3,\ldots, w_{2+L_1}\},\  \{w_{3+L_1},\ldots,
w_{2+\bar{L}_2}\}, \ \ldots, \
\{w_{3+\bar{L}_{k-1}},\ldots,w_{2+\bar{L}_{k}}\}, \ \ldots$$
 as vertices $v^{(L)}_3,v^{(L)}_4,\ldots,v^{(L)}_{k+2},\ldots$.  One interprets
 the sequence of graphs $G_n^{(L)}=\{v_1^{(L)},\ldots,v_{n+2}^{(L)}\}$ for
 $n\geq 0$ as a
more general
 graph process where $L_i$ sets of edges are added
 at the $i$th step preferentially for $i\geq 1$.
This model has some overlap with the very general model given
in \cite{CF}
where vertices can be selected preferentially or
at random; in \cite{CF}, when only ``new'' vertices are selected
preferentially,
their assumptions become $X_i\equiv
1$ and $\{q_j\}_{j\geq 1}$ has bounded support (as well as $\beta =0$).

For the remainder of the article, we will focus, for
simplicity, on the basic ``tree'' model given through
(\ref{weight_eq}), although
extensions to the other case ($L_i\geq 1$, $\beta=0$) under various conditions on
$\{q_j\}_{j\geq 1}$ are possible.

\subsection{Results} For $n\geq 0$ and $j\geq 1$, let
$$R_j(n) \ = \ \sum_{i=1}^{n+2} I\left(d_i(n) = j\right)$$
be the number of vertices in $G_n$ with degree $j$.
Also, define the maximum degree in $G_n$ by
$$
 M_n \ = \ \max_{1\leq i\leq n+2} d_i(n). 
 $$
In addition, denote the mean
$$m \ = \ \sum_{j\geq 1} jp_j.$$

Our first result is on the growth rates of individual degree
sequences $\{d_i(n)\}_{n\geq 0}$ and the maximal one $M_n$. It also describes
the asymptotic behavior of the index where the maximal degree is
attained.

\begin{thm}
\label{Max}
 Suppose $\sum (j\log j)  p_j <\infty$, and let $\theta = m/(2m+\beta)$.

(i) For each $i\geq 1$, there exists a random variable $\gamma_i$ on
$(0,\infty)$ such that
$$
 \lim_{n\rightarrow \infty}
\frac{d_i(n)}{n^\theta} \ = \ \gamma_i \ {\it exists \ a.s.}.
$$

(ii) Further, there exist positive absolutely continuous
independent random variables $\{\xi_i\}_{i\ge 1}$ with $E[\xi_i] <
\infty$, 
and a random variable $V$ on $(0,\infty)$
such that $\gamma_i= \xi_i V$ for $i\geq 1$. In particular, for all
$i,j\geq 1$,

$$\lim_{n \rightarrow \infty}\frac{d_i(n)}{d_j(n)} = \frac{\xi_i}{\xi_j} \;
\textup{\it exists a.s.}.$$

(iii) Also, when $\sum j^rp_j<\infty$ for an $r> \theta^{-1}=
2 +\beta/m$, then
$$\lim_{n\rightarrow \infty} \frac{M_n}{n^\theta}  \ = \
\max_{i\geq 1} \gamma_i \ < \infty\ \ {\it a.s.}
$$

(iv) Moreover, in this case ($\sum j^rp_j<\infty$ for $r>\theta^{-1}$), if $I_n$ is the index where
$$
d_{I_n}(n) \ = \ M_n,
$$
then $\lim_{n\rightarrow \infty} I_n = I <\infty$ exists a.s.
\end{thm}

\begin{remark}
\rm Note that Theorem \ref{Max} asserts that the individual
degrees $d_i(n)$ and the maximal degree $M_n$ grow at the same rate $n^\theta$, and also the vertex with maximal
degree freezes eventually, that is it does not change for large $n$.
\end{remark}

The next result is on the convergence of the empirical distribution
of the degrees $\{d_i(n): 1\le i\le n+2\}$. Let $\{D(y): y\geq 0\}$
be a Markov branching process with exponential($1$) lifetime
distribution, offspring distribution $\{p'_j=p_{j-1}\}_{j\ge 2}$,
immigration rate $\beta\geq 0$, immigration size distribution
$\{p_j\}_{j\geq 1}$, and initial value $D(0)$ distributed according
to $\{p_j\}_{j\ge 1}$
(see Definition \ref{defn-2} in section \ref{sect:embed} for the full
statement). Also, for $y\geq 0$ and $j\geq 1$, let
\begin{equation}
\label{pi_j}  p_j(y) \
= \ P\Big( D(y) = j\Big). \end{equation}

\begin{thm}\label{LLN}
Suppose $\sum (j\log j) p_j<\infty$, and define the probability
$\{\pi_j\}_{j\geq 1}$ by
$$\pi_j \ = \ (2m+\beta) \int_0^\infty p_j(y)e^{-(2m+\beta){y}} dy.$$
Then, for $j\geq 1$, we have
$$
\frac{R_j(n)}{n} \ \rightarrow \ \pi_j, \ \; {\it in
\;probability, \;as } \;  n\rightarrow \infty.
$$
\end{thm}

\begin{remark}\label{cor-bdd}\rm
As a direct consequence, for bounded functions $f: \mathbb{N}
\rightarrow \mathbb{R}$,
$$
\frac{1}{n} \sum_{j=1}^\infty f(j) R_j(n)  \ \rightarrow \
 \sum_{j=1}^\infty f(j)\pi_j, \ \; {\it in \;probability,
\;as } \; n\rightarrow \infty. $$\end{remark}

We now consider the ``power-law'' behavior of the limit degree
 distribution
 $\{\pi_j\}_{j\geq 1}$.
\begin{thm}\label{powerlaw}
Suppose $\sum_{j\geq 1} j^{2+\beta/m} p_j<\infty$.  Then, for $s\geq 0$,
we have
\begin{eqnarray*}
\sum_{j\geq 1} j^s \pi_j < \infty \ \ \ {\rm if \ and \ only \ if \
  \ \ } s<2
+\beta/m.\end{eqnarray*}
\end{thm}

\begin{remark}
\label{power_law_rm}
\rm Heuristically, the last result suggests $\pi_j =
O(j^{-[3+\beta/m]})$ as $j\uparrow \infty$.  In the case
$X_i\equiv x_0$ for $x_0\geq 1$, (\ref{pi_j}) can be
explicitly evaluated (Proposition \ref{explicit_pi}) to
get $\pi_{j} = O(j^{-[3+\beta/x_0]})$ when $j$ is a multiple of $x_0$.
\end{remark}

The next
section discusses the embedding method and auxilliary estimates. In the third
section, the proofs of Theorems \ref{Max}, \ref{LLN}, and
\ref{powerlaw} are given.

\sect{Embedding and some estimates}\label{sect:embed}

We start with the following definitions, and then describe in following
subsections the
embedding and various estimates.

\begin{defn}\label{defn-1}\rm A {\it Markov branching process with offspring
distribution $\{p'_j\}_{j\ge 0}$ and lifetime parameter $0<\lambda<\infty$} is a continuous-time Markov chain
$\{Z(t): t\ge 0\}$ with state space $\S=\{0,1,2\ldots\}$ and waiting
time parameters $\lambda_i\equiv i\lambda$ for $i\geq 0$, and jump probabilities
$p(i,j)= p'_{j-i+1}$ for $j\ge i-1\geq 0$ and $i\geq 1$, $p(0,0)=1$,
and $p(i,j)=0$ otherwise (cf. Chapter III \cite{AN}).\end{defn}

\begin{defn}\label{defn-2}\rm A {\it Markov branching process with offspring
distribution $\{p'_j\}_{j\ge 0}$ and lifetime parameter $0<\lambda<\infty$, immigration parameter $0\leq \beta<\infty$ and immigration size distribution $\{p_j\}_{j\ge
0}$} is a continuous-time Markov chain $\{D(t): t\ge 0\}$ such that
$D(t)=Z(t)$ as in Definition \ref{defn-1} when $\beta=0$, and when $\beta>0$,
$$ D(t)\  =\ \sum_{i=0}^{\infty} Z_{i}(t-T_i)\;I(T_i\leq t)
$$
 where $\{T_i\}_{i\ge 1}$ are the jump times of a Poisson process
$\{N(t): t\ge 0\}$ with parameter $\beta$, $T_0 = 0$, and
$\{Z_i(\cdot)\}_{i\geq 0}$ are independent copies of $\{Z(t):t\ge
0\}$ as in Definition \ref{defn-1}, with $Z_0(0)=D(0)$ and $Z_i(0)$
distributed according to $\{p_j\}_{j\ge 0}$ for $i\geq 1$ and
also independent of $\{N(t): t\ge 0\}$.
\end{defn}

\begin{remark}\label{rem-no-explosion}\emph{The condition that the mean
number of offspring is finite, $\sum j p'_j<\infty$, is sufficient to
ensure that $P(Z(t)< \infty)=1$ and $P(D(t)<\infty)=1$ for all $t\ge 0$, that is no explosion
occurs in finite time (cf. p. 105 \cite{AN})}\end{remark}

\subsection{Embedding process} We now construct a Markov
branching process through which a certain ``embedding'' is
accomplished. Recall $\{p_j\}_{j\geq 1}$ is a probability on the
positive integers. Consider an infinite sequence of independent
processes $\{D_i(t): t\geq 0\}_{i\geq 1}$ where each $\{D_i(t):t\geq
0\}$ is a Markov branching process with immigration
 as in Definition \ref{defn-2}, corresponding to exponential($\lambda
 =1$) lifetimes, offspring distribution $\{p'_j= p_{j-1}\}_{j\geq 2}$
 (with $p'_0=p'_1=0$),
 and immigration parameter $\beta\geq 0$ and immigration size distribution
$\{p_j\}_{j\geq 1}$. The distributions of $\{D_i(0)\}_{i\geq 1}$ will be specified
later.

Now, define recursively the following processes.
\vskip .1cm

$\bullet$ At time $0$, the first
two processes $\{D_i(t): t\geq 0\}_{i=1,2}$ are started with
$D_1(0)=D_2(0)=1$.  Let $\tau_{-1}=\tau_0=0$, and $\tau_1$ be the first time an
``event'' occurs in any one of the two processes.

$\bullet$ Now add a random $X'_1$ of new particles to the process in which the
  event occurred:  (i) If the event is ``immigration,'' then $P(X'_1=j)
  = p_j$ for $j\geq 1$.  (ii) If the event is the death of a particle,
  then $P(X'_1=j)=p_{j-1}$ for $j\geq 2$.  Denote $X_1$ as the net
  addition; then $P(X_1 = j) = p_j$ for $j\geq 1$.

$\bullet$ At time $\tau_1$, start a new Markov branching process $\{D_3(t): t\geq 0\}$ with
  $D_3(0) = X_1$.

$\bullet$ Let $\tau_2$ be the first time after $\tau_1$ that an event occurs
  in any of the processes $\{D_i(t): t\geq \tau_1\}_{i=1,2}$ and
  $\{D_3(t-\tau_1): t\geq \tau_1\}$.  Add a random (net) number $X_2$,
  following the scheme above for $X_1$, of particles with distribution $\{p_j\}_{j\geq 1}$ to the process in
  which the event occurred.  At time $\tau_2$, start a new Markov
  branching process $\{D_4(t): t\geq 0\}$ with $D_4(0) = X_2$.

$\bullet$ Suppose that $n$ processes have been started with the first two at
  $\tau_0=0$, the third at time $\tau_1$, the fourth at time $\tau_2$,
  and so on with the $n$th at time $\tau_{n-2}$, and with (net) additions
  $X_1,X_2,\ldots,X_{n-2}$ at these times.  Now, let $\tau_{n-1}$ be
  the first time after $\tau_{n-2}$ that an event occurs in one of the
  processes $\{D_i(t): t\geq 0\}_{i=1,2}$, $\{D_3(t-\tau_1): t\geq
  \tau_1\}$, $\{D_4(t-\tau_2): t\geq \tau_2\}, \ldots, \{D_n(t-\tau_{n-2}):
  t\geq \tau_{n-2}\}$.  Add a (net) random number $X_{n-1}$ of new
  particles with distribution $\{p_j\}_{j\geq 1}$ (following the scheme above)
to the  process in which the event happened.  Now start the
$(n+1)$st process $\{D_{n+1}(t): t\geq 0\}$ with $D_{n+1}(0) =
X_{n-1}$.

\begin{thm}
\label{prop3_1} {\bf [Embedding Theorem]} Recall the degree sequence
$d_j(n)$ defined for the graphs $\{G_n\}$ near (\ref{weight_eq}). For
$n\ge 0$, let
\begin{eqnarray*} Z_n &\equiv & \{D_j(\tau_n - \tau_{j-2}): 1\leq j\leq
n+2\}, \ \ {\rm and }\\\tilde Z_n &\equiv & \{d_j(n): 1\leq j\leq
n+2\}.\end{eqnarray*} Then, the two collections
$\left\{Z_n\right\}_{n\geq 0}$ and $ \{\tilde Z_n\}_{n\geq 0}$ have
the same distribution.
\end{thm}

{\it Proof.} First note that both sequences $\{Z_n\}_{n\ge 0}$ and
$\{\tilde Z_n\}_{n\ge 0}$ have the Markov  property and
$Z_0=\tilde Z_0=\{1,1\}$. Next, it will be shown below that the
transition probability mechanism from  $Z_n$ to $Z_{n+1}$ is the
same as that from  $\tilde Z_{n}$ to $\tilde
Z_{n+1}$. To see this note that, at time $0$, both $D_1(\cdot)$
and $D_2(\cdot)$ are ``turned on,'' and, at time $\tau_1$,
$D_3(\cdot)$ is ``turned on,'' and more generally, at $\tau_j$,
$D_{j+2}(\cdot)$ is ``turned on.'' At time $\tau_{n+1}$, the
``event'' could be in $D_i(\cdot)$ for $1\leq i\leq n+2$ with
probability
$$\frac{D_i(\tau_n-\tau_{i-2}) +\beta}{\sum_{j=1}^{n+2}
  (D_j(\tau_n -
  \tau_{j-2}) + \beta)}$$
in view of the fact that the minimum of $n+2$ independent
exponential  random variables $\{\eta_i\}_{1\leq i\leq n+2}$ with
means $\{\mu_i^{-1}\}_{1\leq i\leq n+2}$ is an exponential
random variable with mean $(\sum_{i=1}^{n+2}\mu_i)^{-1}$, and
coincides with $\eta_i$ with probability $\mu_i(\sum_{i=1}^{n+2}
\mu_i)^{-1}$ for $1\leq i\leq n+2$. At that event time
  $\tau_{n+1}$, $D_{n+3}(\cdot)$ is ``turned on,'' that is a new
  $(n+3)$rd vertex is created and connected to the chosen vertex
  $v_i$ with $X_{n+1}$ edges between them. Hence both the degree of the new vertex and increment in
  the degree of the chosen vertex (among the existing ones) is
$X_{n+1}$. This shows that the conditional distribution of $Z_{n+1}$
given $Z_{n}=z$ is the same as that of $\tilde Z_{n+1}$ given
$\tilde Z_{n}=z$. \qed \vskip 0.5cm

\subsection{Estimates on branching times}
We now develop some properties of the branching times $\{\tau_n\}_{n\geq 1}$, used in the
embedding in subsection 2.1, which have some analogy to results in section III.9
\cite{AN} (cf. \cite{athreya-karlin}).
Define $S_0 = 2 +2\beta$ and, for $n\geq 1$,
 $$S_n \ = \ 2+2\beta + \sum_{j=1}^n 2X_j + n\beta,$$ where as before $X_1,\ldots,
X_n$ are the independent and identically distributed according to
$\{p_j\}_{j\geq 1}$ net additions at event times $\tau_1,\ldots,\tau_n$.

\begin{prop}
\label{prop1} The random variable $\tau_1$ is exponential with
mean $S_0^{-1}$.  Also, for $n\geq 1$, conditioned on the
$\sigma$-algebra ${\cal F}_n$ generated by $\{D_j(t-\tau_{j-2}):
\tau_{j-2}\leq t\leq \tau_n; X_j\}_{1\leq j\leq
  n}$, the random variable $\tau_{n+1}-\tau_n$
is exponential with mean $S_n^{-1}$.
\end{prop}

{\it Proof.} Follows from the construction of the $\{\tau_i\}_{i\geq 1}$. \qed
\vskip 0.5cm

\begin{prop}
\label{prop2} Suppose $m = \sum jp_j<\infty$.
Then,
$$\bigg\{\tau_n - \sum_{j=1}^n \frac{1}{S_{j-1}}; {\cal
  F}_n\bigg\}_{n\geq 1}$$
is an $L^2$ bounded martingale and hence converges a.s. as well as
in $L^2$.
\end{prop}

{\it Proof.}
The martingale property follows from the fact
$$\tau_n \ = \ \sum_{j=1}^n \big( \tau_j - \tau_{j-1} \big)$$
and Proposition \ref{prop1}.

Next, with $\phi(a) = E[e^{-aX_1}]$ for $a\geq 0$, we have the uniform
bound in $n\geq 1$,
\begin{eqnarray*}
{\rm Var}\bigg( \tau_n - \sum_{j=1}^n \frac{1}{S_{j-1}}\bigg) & =&
{\rm Var}\bigg( \sum_{j=1}^n \left(   \tau_j - \tau_{j-1}
- \frac{1}{S_{j-1}}\right)\bigg)\\
& =& \sum_{j=1}^n {\rm Var}\bigg( \tau_j - \tau_{j-1} -
\frac{1}{S_{j-1}}\bigg)\;\;\;\;\;\textup{ (by martingale property)}\\
& =& \sum_{j=1}^n E\bigg[ \frac{1}{S^2_{j-1}}\bigg]\\
& =& \sum_{j=1}^n E\bigg[ \int_0^\infty x e^{-S_{j-1} x}dx\bigg]\\
& \leq& \sum_{j=1}^\infty \int_0^\infty \bigg(\phi(2x)
e^{-x\beta}\bigg)^{j-1} xe^{-(2 +2\beta)x}dx \\
& \leq& \int_0^\infty \frac{x
  e^{-(2 +2\beta)x}dx}{1-\phi(2x)e^{-x\beta}} \ \
<\infty
\end{eqnarray*}
where the finiteness in the last bound follows from the fact that
\[\lim_{x\downarrow 0}\frac{x}{1-\phi(2x)e^{-x\beta}} =
\frac{1}{2m+\beta} \
<\infty.\] The a.s. and $L^2$-convergence follows from Doob's
martingale convergence theorem (c.f. Theorem 13.3.9 \cite{AL}).\qed \vskip 0.5cm

\begin{prop}\label{prop3}
Suppose $\sum (j\log j) p_j<\infty$, and recall $m=\sum jp_j$.
Let also $\alpha = (2m +\beta)^{-1}$.  Then, there exists a
real random variable $Y$ so that a.s.,
$$\lim_{n\rightarrow \infty} \tau_n - \sum_{j=1}^n
\frac{\alpha}{j} \ = \ Y.$$
\end{prop}

{\it Proof.}  By Proposition \ref{prop2}, there is a
real random variable $Y'$ such that,
$$\tau_n - \sum_{j=1}^n \frac{1}{S_{j-1}}\ \rightarrow \ Y' \ \ \ {\rm
  a.s.}
$$
To complete the proof, we note, as $E[X_1\log X_1]= \sum (j\log j)p_j<\infty$,
by Theorem III.9.4 \cite{AN} on reciprocal sums, that
$\sum_{j=1}^\infty ({1}/{S_j} - {\alpha}/{j})$ converges a.s. \qed


\begin{cor}\label{new-cor} Suppose $m=\sum jp_j<\infty$.  Then,

\noindent (i)  $\tau_n \uparrow \infty  \; \; \textup{\it a.s.,
 as } n \rightarrow \infty.$

Also, when $\sum (j\log j)p_j<\infty$, we have, with $\alpha =
(2m+\beta)^{-1}$, that

\noindent (ii) $\tau_n - \alpha\log n \rightarrow \tilde Y := Y - \alpha
\gamma  \textup{ \it a.s.,   as } n \rightarrow \infty,$ where
$\gamma$ is the Euler's
constant.\\
(iii) For each fixed $\epsilon >0$, $\; \sup_{n\epsilon\le k\le
n}\left(\tau_n - \tau_k - \alpha\log (n/k) \right)\rightarrow 0
\textup{ \it a.s.,   as } n \rightarrow \infty$.
\end{cor}

{\it Proof.} The first claim follows from Proposition
\ref{prop2} and the fact that $\sum 1/S_j = \infty$, since by strong
law of large numbers, we have a.s. that $S_j \le j(1/\alpha +1)$ for
large $j$.  The last two claims, as $\sum_{j=1}^n 1/j - \log n
\rightarrow \gamma$, Euler's constant, are direct consequences
of Proposition \ref{prop3}.\qed \vskip 0.5cm

\subsection{Estimates on Markov branching processes}

As in Definition \ref{defn-2}, let $\{D(t): t\geq 0\}$ be a Markov
branching process with offspring distribution $\{p'_j =
p_{j-1}\}_{j\geq 2}$, lifetime $\lambda =1$ and immigration $\beta\geq
0$ parameters, and immigration distribution $\{p_j\}_{j\geq 1}$.

\begin{prop}
\label{prop5}  
Suppose $\sum (j\log j) p_j<\infty$, and $D(0)\geq 1$,
$E[D(0)]<\infty$.  Recall $m=\sum jp_j$.  Then,
$$\lim_{t\rightarrow \infty} D(t)e^{-mt} \ = \ \zeta$$
converges a.s. and in $L^1$, and $\zeta$ is supported on $(0,\infty)$
and has an absolutely continuous
distribution.
\end{prop}

{\it Proof.} Let $\beta>0$; when $\beta =0$ the argument is easier and
a special case of the following development.  Let $0=T_0<T_1<\cdots T_n< \cdots$ be the times at
which immigration occurs, and let $\eta_1,\eta_2, \ldots$ be the
respective number of immigrating individuals (distributed according
to $\{p_j\}_{j\geq 1}$).  From Definition \ref{defn-2}, $D(t)$ has
representation
\begin{equation}
\label{prop5_1} D(t) \ = \ \sum_{i=0}^\infty Z_i\big(t- T_i\big)\;
I\big(T_i\leq t\big)\end{equation} where $\{Z_i(t): t\geq 0\}_{i\geq
0}$ are independent Markov branching processes with offspring
distribution $\{p'_j = p_{j-1}\}_{j\geq 2}$, with exponential($\lambda
=1$) lifetime distributions, with no immigration, with
$Z_0(0) = D(0)$ and $Z_i(0)=\eta_i$ for $i\geq 1$, and also
independent of $\{T_i\}_{i\geq 0}$.
Under the hypothesis $\sum (j\log j)p_j<\infty$, it is known
(Theorem III.7.2 \cite{AN}; with rate $\lambda(\sum_{j\geq 2} jp'_j -1) =
\sum_{j\geq 1}(j+1)p_j -1 = m$), for $i\geq 0$, that
\begin{equation}
\label{prop5_2}\lim_{t\rightarrow \infty} Z_i(t) e^{-mt} \ =\  W_i\end{equation}
  converges in $(0,\infty)$ a.s. and $W_i$ has a continuous
  distribution on $(0,\infty)$. Also  under the hypothesis that $\sum
  (j \log j) p_j < \infty$, it can be shown (Proposition \ref{appendix})
 that
\begin{equation}
\label{prop5_3}
E[\widetilde{W}_i]<\infty \ \  {\rm where \ \ }\widetilde{W}_i \ = \ \sup_{t\geq 0} Z_i(t)
e^{-mt},\end{equation}
and hence convergence in (\ref{prop5_2}) holds in $L^1$ as well.

Since $\{T_i\}_{i\geq 0}$ is a Poisson process with rate $\beta$,
and independent of $\{Z_i(t)\}_{t\ge 0}$,
\begin{eqnarray}\label{prop5_4}E\bigg[\sum_{i=0}^\infty
  \widetilde{W}_i e^{-mT_i}\bigg] \
\leq \ E[\widetilde{W}_1]\bigg( E [D(0)] + \sum_{i=1}^\infty
\bigg(\frac{\beta}{m+\beta}\bigg)^i \bigg)\ < \
\infty,\end{eqnarray} yielding
\begin{equation}
\label{tilde(W)_sum_convergence}
\sum_{i=0}^\infty \widetilde{W}_i e^{-mT_i} \ < \ \infty \ \ \ {\rm
  a.s.}.\end{equation}
Hence, noting (\ref{prop5_2}), (\ref{prop5_3}) and (\ref{tilde(W)_sum_convergence}), by dominated convergence,
\begin{eqnarray}\label{prop5_5}
\lim_{t\rightarrow \infty} D(t)e^{-mt} &=& \sum_{i=0}^\infty
\lim_{t\rightarrow \infty} \bigg[Z_i(t-T_i)I(T_i\leq
t)e^{-m(t-T_i)}\bigg]e^{-mT_i}\nonumber \\
&=& \sum_{i=0}^\infty W_i e^{-mT_i} \ := \ \zeta
\end{eqnarray}
converges in $(0,\infty)$ a.s.. Also,
\begin{equation}\label{prop5_tilde_D}
\sup_{t\ge 0} D(t) e^{-mt}
\ \le \ \sum_{i=0}^\infty \widetilde{W}_i e^{-mT_i}\end{equation}
and hence by
(\ref{prop5_4}) and (\ref{prop5_5}), we get that
\begin{eqnarray*}
\lim_{t\rightarrow \infty} D(t)e^{-mt} &=&  \zeta \;\;\textup{ in
$L^1$.}\end{eqnarray*}

Finally, since
$\{W_i\}_{i\ge 0}, \{T_i\}_{i\ge 1}$ are
independent, absolutely continuous random variables, $\zeta$ is absolutely continuous as well. \qed \vskip
0.5cm

\subsection{Suprema estimates}

We give now some moment estimates which follow by
combination of results in the literature.
Let $\{Z(t):t\geq 0\}$ be a Markov branching process with offspring
distribution $\{p'_j=p_{j-1}\}_{j\geq 2}$ and lifetime parameter $\lambda=1$
as in Definition \ref{defn-1} with independent initial population $Z(0)$
distributed according to $\{p_j\}_{j\geq 1}$. Recall $m=\sum jp_j$,
and, from (\ref{prop5_2}) and (\ref{prop5_3}), that
$$W\ = \ \lim_{t\rightarrow \infty}Z(t) e^{-mt} \ \ {\rm and \ \ }
\widetilde{W} \ = \ \sup_{t\geq 0} Z(t)e^{-mt}.$$
\begin{prop}
\label{appendix}
The following implications hold:
$${\rm If \ }\sum_{j\geq 1} (j\log j)p_j<\infty \ {\rm \ \ then \ \ \ } E[\widetilde{W}]<\infty.$$
Also, for $s>1$,
$${\rm if}  \ \sum_{j\geq 1} j^s p_j<\infty \ \ \ {\rm then \ \ \ } E[\widetilde{W}^s]<\infty.$$
\end{prop}

{\it Proof.} We will take without loss of generality $Z(0)=1$, as
the initial value $Z(0)$ in both statements of the proposition is assumed
to have enough integrability.  Then, first, as $Z(\cdot)$ is
increasing, we have
\begin{equation}
\label{app1}\widetilde{W} \ = \ \sup_{t\geq 0} Z(t)e^{-mt} \ \leq \
e^{m}\sup_{n\geq
  0} Z(n)e^{-mn}\ := \ \widetilde{W}_0.\end{equation}

The process $\{Z(n): n\geq 0\}$ is a discrete-time branching, $W=
\lim_{n\rightarrow \infty} Z(n)e^{-mn}$, and $P(Z(1)=0)=0$,
$P(Z(1)=j)<1$ for all $j\geq 1$. From Lemma I.2.6 in Asmussen
\cite{Asmussen}, for $r\geq 1$, we have, when $P(W>0)>0$, that
\begin{equation}
\label{app2}E[\widetilde{W}_0^r] \ \leq  \ C_0\bigg(1+ E[W^r]\bigg)\end{equation}
for a constant $C_0$.

From Theorem I.10.1 \cite{AN} or Theorem I.2.1 \cite{Asmussen},
\begin{equation}
\label{app3}
P(W>0)>0, \ E[W]<\infty \ \ \ {\rm if \ and \ only \ if \ }  \ \ \sum_{j\geq 1}
(j\log j)P(Z(1)=j) <\infty.\end{equation}
In particular, when $\sum j^rP(Z(1)=j)<\infty$ for $r>1$, $P(W>0)>0$.

Also,
from Theorem I.4.4 \cite{Asmussen}, and the discussion on p. 41-42 \cite{Asmussen}
(cf. equation (4.15) \cite{Asmussen}), we have, for $r>1$ 
and two constants $C_1,C_2$,
that
\begin{equation}\label{app4} C_1 \sum_{j\geq 1} j^r P(Z(1)=j)\ \leq \ E[W^r]  \ \leq
C_2\bigg(1+ \sum_{j\geq 1} j^r P(Z(1)=j)\bigg);\end{equation} hence,
$E[W^r]<\infty$ when $\sum j^r P(Z(1)=j)<\infty$.

From Corollary III.6.2 \cite{AN} (cf. \cite{athreya}), for
$a\geq 1$, $b\geq 0$,
\begin{equation}\label{app5}\sum_{j\geq 1} j^a|\log j|^b P(Z(1)=j)
  <\infty
\ \ \ {\rm if \ and \ only \ if \ }
\ \  \sum_{j\geq 1} j^a|\log j|^b p_j
<\infty.\end{equation}

Then, straightforwardly combining (\ref{app1})-(\ref{app5}), we
conclude the proof. \qed

Let now $\{D(t):t\geq 0\}$ be a Markov branching process with
offspring distribution $\{p'_j=p_{j-1}\}_{j\geq 2}$, lifetime $\lambda=1
$ and immigration $\beta\geq 0$ parameters, and immigration distribution
$\{p_j\}_{j\geq 1}$ as in Proposition \ref{prop5} with
also $D(0)$ distributed as $\{p_j\}_{j\geq 1}$.  Let also
$$\widetilde{D} \ := \ \sup_{t\geq 0} D(t)e^{-mt}.$$

\begin{prop}
\label{appendix1}
For $r>1$, we have
$${\rm if \ }\sum_{j\geq 1} j^r p_j <\infty \ {\rm \ \ then } \ \ \ E[\widetilde{D}^r] <
 \infty.$$
\end{prop}

{\it Proof.}  When $\beta =0$, the statement is the same as
Proposition \ref{appendix}.  When $\beta>0$, as in the proof of Proposition \ref{prop5}, let
$\{T_i\}_{i\geq 1}$ be the times of immigration, and $T_0=0$.  Note that
$\sum_{i\geq 0} e^{-mT_i} <\infty$ a.s. as the expected value
$\sum_{i\geq 0} (\beta/(m+\beta))^i$ is
finite.  From (\ref{prop5_tilde_D}), and Jensen's inequality, we have
\begin{eqnarray*}
\widetilde{D}^r &\leq & \bigg(\sum_{i\geq 0}
\widetilde{W}_i
  e^{-mT_i}\bigg)^r\\
&\leq & \bigg(\bigg[\sum_{j\geq 0} e^{-mT_j}\bigg]^{-1} \sum_{i\geq 0}
  \widetilde{W}_i^r e^{-mT_i}\bigg)\bigg( \sum_{j\geq 0}
  e^{-mT_j}\bigg)^r\\
&=& \bigg( \sum_{i\geq 0}\widetilde{W}_i^r
e^{-mT_i}\bigg)\bigg(\sum_{j\geq
  0} e^{-mT_j}\bigg)^{r-1}.\end{eqnarray*}
Hence, by independence of $\{\widetilde{W_i}\}_{i\geq 0}$ and
  $\{T_i\}_{i\geq 0}$, for an integer $K\geq r-1$, we have
$$E\big[\widetilde{D}^r\big] \ \leq \  E[\widetilde{W}^r_1]\sum_{i\geq
  0}E\bigg[ e^{-mT_i}\bigg(\sum_{j\geq
  0} e^{-mT_j}\bigg)^{K}\bigg].$$
From Proposition \ref{appendix}, $E[\widetilde{W}^r_1]<\infty$.  Also,
\begin{eqnarray*}E\bigg[ e^{-mT_i}\bigg(\sum_{j\geq
  0} e^{-mT_j}\bigg)^{K}\bigg] & \leq & E[e^{-2mT_i}]^{1/2}E\bigg[\bigg(\sum_{j\geq
  0} e^{-mT_j}\bigg)^{2K}\bigg]^{1/2}\\
&=& \bigg(\sqrt{\frac{\beta}{2m+\beta}}\bigg)^{i}E\bigg[\bigg(\sum_{j\geq
  0} e^{-mT_j}\bigg)^{2K}\bigg]^{1/2}.\end{eqnarray*}
To finish, we now bound, given $T_j$ is the sum of $j$ independent exponential
  random variables with parameter $\beta$ for $j\geq 1$, that
\begin{eqnarray*}
&&E\bigg[\bigg(\sum_{j\geq
  0} e^{-mT_j}\bigg)^{2K}\bigg]\\
 &&\ \ \ \ \ \ = \ (2K)! \sum_{0\leq j_1\leq
  \cdots \leq \ j_{2K}} E\bigg[ \prod_{l=1}^{2K} e^{-mT_{j_l}}\bigg]\\
&&\ \ \ \ \ \ = \ (2K)!\sum_{0\leq j_1\leq \cdots \leq j_{2K-1}}
E\bigg[ \prod_{l=1}^{2K-2} e^{-mT_{j_l}}
  e^{-2mT_{j_{2K-1}}}\bigg] \bigg
  (\frac{\beta}{m+\beta}\bigg)^{j_{2K}-j_{2K-1}}\\
&&\ \ \ \ \ \ \leq \  (2K)! \bigg(\frac{m+\beta}{m}\bigg)\sum_{0\leq j_1\leq \cdots \leq j_{2K-1}}
E\bigg[ \prod_{l=1}^{2K-1} e^{-mT_{j_l}}\bigg] \\
&&\ \ \ \ \ \ \leq \ (2K)! \bigg(\frac{m+\beta}{m}\bigg)^{2K}\end{eqnarray*}
is finite for fixed $K$. \qed

\sect{Proof of main results 
}

We give the proofs of the three main results in successive subsections.

\subsection{
Growth rates for degrees and the maximal degree}

We first begin with a basic analysis result.

\begin{prop}
\label{prop_max} Let $\{a_{n,i}: 1\leq i\leq n\}_{n\geq 1}$ be a
double array of nonnegative numbers such that

(1) For all $i\geq 1$, $\lim_{n\rightarrow \infty} a_{n,i} = a_i<\infty$,

(2) $\sup_{n\geq 1} a_{n,i} \leq b_i<\infty$  and

(3) $\lim_{i\rightarrow \infty} b_i = 0$. 


\noindent Then,

(a) $\max_{1\leq i\leq n} a_{n,i} \rightarrow \max_{i\geq 1} a_i$, as
$n\rightarrow \infty$. 

(b) In addition, if $a_i \neq a_j$ for distinct $i,j\geq 1$, there exists $I_0$ and $N_0$ such that $\max_{1\leq i\leq n} a_{n,i} = a_{n,I_0}$
for $n\geq N_0$.
\end{prop}

{\it Proof.}  For each $k\geq 1$,
$$\lim_{n\rightarrow \infty}\max_{1\leq i\leq k}a_{n,i} \ = \ \max_{1\leq i\leq k} a_i.$$
Hence,
$$
\liminf_{n\rightarrow \infty} \max_{i\geq 1} a_{n,i}  \ \geq \
\liminf_{n\rightarrow \infty} \max_{1\leq i\leq
  k} a_{n,i} \ = \ \max_{1\leq i\leq k} a_i$$ which gives
\begin{equation}
\label{prop7_1}
\liminf_{n\rightarrow \infty} \max_{i\geq 1} a_{n,i}  \ \geq \  \max_{i\geq 1} a_i. \end{equation}

Also, for each $k\geq 1$,
$$\max_{i\geq 1} a_{n,i} \ \leq \ \max_{1\leq i\leq k} a_{n,i} +
\max_{i>k}b_i.$$
Then,
$$\limsup_{n\rightarrow \infty}\max_{i\geq 1} a_{n,i} \ \leq \ \max_{1\leq i\leq k} a_i +
\max_{i>k} b_i \ \leq \ \max_{i\geq 1} a_i + \max_{i\geq k} b_i.$$
Since $\lim_{i\rightarrow \infty} b_i = \limsup_{i\rightarrow \infty}
b_i =
\lim_{k\rightarrow \infty} \max_{i\geq k}b_i =0$, we have
\begin{equation}
\label{prop7_2}
\limsup_{n\rightarrow \infty} \max_{i\geq 1} a_{n,i} \ \leq \ \max_{i\geq 1} a_i.
\end{equation}

Now, (\ref{prop7_1}) and (\ref{prop7_2}) yield part (a).  By
assumption (3), $\max_{i\geq 1} a_i$ is attained at some finite index $I_0$, and
by assumption (4) this index is unique, giving part (b).\qed
\vskip 0.5cm

\noindent{\bf Proof of Theorem \ref{Max}.} By the embedding theorem
(Theorem \ref{prop3_1}), to establish Theorem \ref{Max} for the
sequence $\{\widetilde Z_n\}_{n\ge 0}$, it suffices to prove the
corresponding results for the $\{Z_n\}_{n\ge 0}$ sequence.

By Proposition \ref{prop5} and Corollary
\ref{new-cor}(i),
$$\lim_{n\rightarrow \infty}
D_i(\tau_n-\tau_{i-2})e^{-m(\tau_n-\tau_{i-2})} \ = \ {\zeta}_i$$
converges a.s. in $(0,\infty)$ for $i\geq 1$.  By Proposition \ref{prop3}, a.s. as $n\uparrow
\infty$,
\begin{eqnarray*}
\exp\bigg\{\frac{-m}{2m+\beta}\sum_{j=1}^n
\frac{1}{j}\bigg\}\exp\big\{ m\tau_n\big\} &=&
\exp\bigg\{m\bigg(\tau_n-\sum_{j=1}^n
\frac{1}{j(2m+\beta)}\bigg)\bigg\}\\
&\rightarrow & e^{mY}.\end{eqnarray*} Further, $\sum_{j=1}^n(1/j) -
\log n \rightarrow \gamma$, Euler's constant.  Thus, a.s. as
$n\uparrow \infty$,
\begin{equation}
\label{tau_n_conv}
e^{m\tau_n} n^{-m/(2m+\beta)} \ \rightarrow \
e^{mY}e^{m\gamma/(2m+\beta)}\ := \ V\end{equation}
where $V$ is a positive real random variable.
 Hence,
\begin{eqnarray*}
D_i(\tau_n-\tau_{i-2})n^{-m/(2m+\beta)} &=&
D_i(\tau_n-\tau_{i-2})e^{-m(\tau_n-\tau_{i-2})} e^{-m\tau_{i-2}}\\
&&\ \ \ \ \ \ \ \ {} \hskip 3cm \times
e^{m\tau_n}n^{-m/(2m+\beta)}\\
&\rightarrow& \ {\zeta}_i e^{-m\tau_{i-2}}V \;:=\; \xi_i V,
\end{eqnarray*} a.s.  as $n\uparrow \infty$,  where
$\xi_i={\zeta}_ie^{-m\tau_{i-2}}$ is a positive real random variable.
This proves part (i) with $\gamma_i = \xi_i V$.

By independence of $\tau_{i-2}$ and $\{D_i(t)\}_{t\ge 0}$, absolute
continuity of $\tau_{i-2}$ for $i\geq 3$ ($\tau_0=\tau_{-1}=0$), and Proposition
\ref{prop5},
it
follows that $\xi_i$ has an absolutely continuous distribution with
finite mean,  proving part (ii).

To prove part (iii) and (iv), we
first note, for each $i\geq 1$, that
$$D_i(\tau_n-\tau_{i-2})e^{-m(\tau_n-\tau_{i-2})}
 \ \leq \ \sup_{t\geq 0} D_i(t)e^{-mt}\ := \ \widetilde{D}_i.$$
Let
\begin{eqnarray*}
a_{n,i} &=& D_i(\tau_n - \tau_{i-2})e^{-m\tau_n} \ \ {\rm for \ } 1\leq
i\leq n, \ \ {\rm and} \\
b_i&=& \widetilde{D}_i e^{-m\tau_{i-2}} \ \ {\rm for \ } i\geq
1.\end{eqnarray*} For each $i\geq 1$, $\sup_{n\geq 1}a_{n,i} \leq
b_i$ and $a_{n,i} \rightarrow {\zeta}_ie^{-m\tau_{i-2}}:=a_i$ say.
Since $\sum_{j=1}^{\infty} j^r p_j < \infty$ for some $r > 1$
(satisfying $rm/(2m+\beta)=r\theta >1$), we have that $E(\widetilde{D}_i^r)< \infty$
(Proposition \ref{appendix1}).
 By Markov's inequality, for
all $\epsilon >0$,
%
$$ P\left(\widetilde{D}_i > \epsilon
i^{m/(2m+\beta)}\right) \ \leq\
E[\widetilde{D}_1^{r}]/ \big(\epsilon^r
i^{rm/(2m+\beta)}\big).$$
Hence, by Borel-Cantelli, we have a.s.
$$\widetilde{D}_i \le \epsilon i^{m/(2m+\beta)} \textup{ for all large }
i.$$ Since, by Corollary \ref{new-cor}(ii), $m\tau_k$ is on order
$[{m}/{(2m+\beta)}]\log k$ a.s. for large $k$, and $\epsilon>0$ is
arbitrary, it follows that $b_i = \widetilde{D}_i
e^{-m\tau_{i-2}} \rightarrow 0$ a.s. as $i\uparrow \infty$. By
Proposition \ref{prop_max} and (\ref{tau_n_conv}),
this implies that a.s.,
$$\lim_{n\rightarrow \infty}
\max_{1\leq i\leq n}D_i(\tau_n-\tau_{i-2})n^{-m/(2m+\beta)}
 \ =\
V\max_{i\geq 1} \zeta_i e^{-m\tau_{i-2}}.$$

Now we claim $\{\zeta_i e^{-m\tau_{i-2}}\}_{i\geq 1}$ are all
distinct, that is $P(\zeta_j e^{-m\tau_{j-2}}=\zeta_i
e^{-m\tau_{i-2}})= P(\zeta_j e^{-m(\tau_{j-2}-\tau_{i-2})}=\zeta_i)=
0$ for any $1\leq i < j$. This follows from the fact that
conditioned on ${\cal F}_{i-2}$ (see Proposition \ref{prop1} for
definition of ${\cal F}_{i}$) the random variables, $\zeta_j$ and
$\tau_{j-2}-\tau_{i-2}$, are independent with
absolutely continuous distributions (when $j\geq 3$ for
$\tau_{j-2}-\tau_{i-2}$ as $\tau_{0}=\tau_{-1}=0$) and are independent of
$\zeta_i$.
Hence, $V\max_{i\geq 1} \zeta_i e^{-m\tau_{i-2}}$ is attained at a
unique index $I_0$. Also, as $\{I_n\}_{n\geq 1}$ are
  integer valued random variables, $I_n$ will equal $I_0$ a.s. for all
  large $n$.
\qed
%



\subsection{ Convergence of the empirical distribution of degrees} 


The following lemma will be helpful in the proof of Theorem
\ref{LLN}.


\begin{lem}\label{new-lemma} Let $\{X(t): t\ge 0\}$ be
a continuous-time, discrete state-space, Markov chain
which is non-explosive, that is the number of jumps of
$\{X(t):t\ge 0\}$  in any finite time-interval $[0,K]$ is finite
a.s.  For $K>0,\; \delta>0$, let
\[p_K(\delta)\;\equiv \;\sup_{0\le t\le K} P\Big( |
X(t+\delta)-X((t-\delta)\vee 0)| \ge 1\Big).\] Then, for all $ K
>0$,
\[\lim_{\delta \downarrow 0} \ p_K(\delta)=0.\\\]
\end{lem}

{\it Proof.} Since $\{X(t): t\ge 0\}$  is non-explosive, for any
$0<K<\infty$, the number of jumps $N(K)$ of $\{X(t): 0\le t\le K\}$
is a finite valued random variable a.s.. Also for any $j < \infty$,
the jump times $(T_1, \ldots, T_j)$ of $\{X(t): t\ge 0\}$  have a
continuous joint distribution. These two facts together yield the
lemma. \qed

The following result follows from Remark \ref{rem-no-explosion} and the above lemma.

\begin{cor}
\label{new-lemma-cor} Let $\{D(t): t\geq 0\}$ be as
in Definition \ref{defn-2}.
Define, for $0\leq s\leq t$,
$D(s,t]=D(t)-D(s)$ and, for $K>0, \delta>0$,
\[p_K^{D}(\delta)\;= \;\sup_{0\le t\le K} P\Big(D((t-\delta)\vee 0, t+\delta]
\ge 1\Big).\] Then for $K >0$,
\[\lim_{\delta \downarrow 0} \ p_K^D(\delta)=0.\\\]
\end{cor}


\noindent{\bf Proof of Theorem \ref{LLN}.}  Recall $\alpha = (2m +\beta)^{-1}$ (from subsection 2.2).  For $n\geq 1$,
note that \begin{eqnarray}\frac{R_j(n)}{n} &= &
\frac{1}{n}\sum_{i=1}^{n+2} I\big(D_i(\tau_n - \tau_{i-2}) =
j\big)\nonumber\\
{}&=& \frac{1}{n}\sum_{i=1}^{n+2} \Big\{I\left(D_i(\tau_n -
\tau_{i-2}) = j\right) -I\left(D_i\left(\alpha
\log\left({n}/{(i-2)}\right)\right) = j\right)\Big\}\label{a1}\\
{}&&+ \;\frac{1}{n}\sum_{i=1}^{n+2} \Big\{I\left(D_i\left(\alpha
\log\left({n}/{(i-2)}\right)\right) = j\right) -p_j\left(\alpha
\log({n}/{(i-2)})\right)\Big\}\label{a2}\\
{}&&+ \;\frac{1}{n}\sum_{i=1}^{n+2} \Big\{p_j\left(\alpha
\log({n}/{(i-2)})\right) -\int_0^1 p_j\left(-\alpha
\log(x)\right)\Big\}\;\label{a3}\\
{}& &+\; \int_0^1 p_j\left(-\alpha \log(x)\right) dx\label{a4}\\
{}&=& J_1(n) + J_2(n) + J_3(n) +J_4(n)\nonumber,
\end{eqnarray}
 where $J_i(n), i=1,\ldots,4$ represent the $4$ terms in
(\ref{a1}-\ref{a4}) (for notational convenience, we use the
convention $\log(n/(i-2))\equiv \log(n)$ for $i=1,2$ here). The
proof is now obtained by showing $J_i(n)$ vanishes in probability for
$i=1,2,3$ and observing, after change of variables, that
$$J_4(n) \ \equiv \ \frac{1}{\alpha}\int_0^\infty p_j(y)
\exp\{-y/\alpha\}dy.$$

To show that the first term $J_1(n)$ goes to 0 in probability, fix
$\epsilon
>0$. Note, from Corollary \ref{new-cor}(iii), for $\delta >0$, if
$${\cal A}_n(\delta) \equiv \left\{\sup_{
n\epsilon+2 \le i\le n+2} \big | \tau_n - \tau_{i-2} - \alpha
\log(n/(i-2))
\big | > \delta \right\},$$
then,
\begin{equation}
\label{a-n-lln}
\limsup_{\delta \downarrow 0}\limsup_{n\rightarrow \infty} P\left({\cal
A}_n(\delta)\right)\ =\ 0.\end{equation}
\noindent Now for $n\epsilon +2\leq i\leq n+2$. we have $0 \leq
\alpha \log\left({n}/{(i-2)}\right) \le -\alpha \log \epsilon $.
Hence, from the definition of ${\cal A}_n(\delta)$, we get the
following bound on the expectation of a typical summand in $J_1$ in
this range.
%
\begin{eqnarray*}
&&
E\left( \big | I\left(D_i(\tau_n - \tau_{i-2}) = j\right)
-I\left(D_i\left(\alpha \log\left({n}/{(i-2)}\right)\right) =
j\right)\big |\right) \\
&&\ \ \ \ =\  P\left( \big | I\left(D_i(\tau_n - \tau_{i-2}) = j\right)
-I\left(D_i\left(\alpha \log\left({n}/{(i-2)}\right)\right) =
j\right)\big
|=1\right)\\
&&\ \ \ \ \leq \ P\left( \big | D_i(\tau_n - \tau_{i-2}) - D_i\left(\alpha
\log\left({n}/{(i-2)}\right)\big
| \ge 1\right)\right)\\
&&\ \ \ \ \leq \ P \Big( \{ | D_i(\tau_n - \tau_{i-2}) - D_i\left(\alpha
\log\left({n}/{(i-2)}\right)\right) | \ge 1 \}
\cap {\cal A}_n^c(\delta)\Big) + P\left({\cal A}_n(\delta)\right)\\
&&\ \ \ \ \le  \ P \Big(
\{D_i\big((\alpha\log\left({n}/{(i-2)}\right)-\delta )\vee 0,\alpha\log\left({n}/{(i-2)}\right)+
\delta\big] \ge 1\}\cap {\cal A}_n^c(\delta)\Big)\\
&&\ \ \ \ \ \ \ \ \ \ \ \ \ \ \  + P\left({\cal A}_n(\delta)\right)\\
&&\ \ \ \ \le\   \sup_{a\in [0,-\alpha\log \epsilon]}P \Big(
\{D_i\big((a-\delta)\vee 0,a+ \delta\big] \ge 1\}\cap {\cal
A}_n^c(\delta)\Big)
+ P\left({\cal A}_n(\delta)\right)\\
&&\ \ \ \ = \ p^{D_3}_{K(\epsilon)}(\delta) + P\left({\cal A}_n(\delta)\right)
\end{eqnarray*}
where $K(\epsilon)=-\alpha\log \epsilon$, and we recall $D_3(\cdot)$ is a
Markov branching process with immigration with $D_3(0)$ distributed
according to $\{p_j\}_{j\geq 1}$. Since each summand in
$J_1$  is bounded  (by $1$),  we have by splitting the sum over
indices $1\leq i<n\epsilon +2$ and $n\epsilon +2\leq i\leq n+2$, the
following bound:
\begin{eqnarray*} E\big|J_1(n)\big| \le \frac{1}{n}(n\epsilon+2) +
\frac{1}{n} (n-n\epsilon)\left[p^{D_3}_{K(\epsilon)}(\delta) +
P\left({\cal A}_n(\delta)\right)\right]
\end{eqnarray*}
Now for fixed $\epsilon >0$, taking limit as $n \rightarrow \infty$
first and then over $\delta \rightarrow 0$, we get from Corollary
\ref{new-lemma-cor} and (\ref{a-n-lln}) that
\[\limsup_{n\rightarrow\infty} E|J_1(n)| < \epsilon\]
and, as $\epsilon >0$ is arbitrary, that $\lim_{n\rightarrow
  \infty} E|J_1(n)|=0$. Hence $J_1(n) \rightarrow  0$
in probability, as $n\rightarrow \infty$.

For the second term $J_2(n)$, we have from Markov's inequality that
for any $\epsilon >0$,
\begin{eqnarray*}
&&P\left( |J_2(n)|
>\epsilon\right) \\
&&\ \ \ \ \ \leq \  \frac{1}{n^4\epsilon^4}
E\bigg(\sum_{i=1}^{n+2} \big(I\big(D_i(\alpha
\log\left({n}/{(i-2)}\right))=j\big) - p_j\big(\alpha
\log\left({n}/{(i-2)}\right)\big) \bigg)^4\\
&&\ \ \ \ \ \leq\  \frac{6}{n^2
\epsilon^4}, \end{eqnarray*} using independence of
$\{D_i(\cdot)\}_{i\geq 1} $, and hence of  the summands above.  Now,
by Borel-Cantelli arguments and the
 method of fourth moments (cf. Theorem 8.2.1 \cite{AL}), we get
 $J_2(n)\rightarrow 0$ a.s., as $n\rightarrow \infty$.

Finally, by simple estimates, and Riemann integrability of
$p_j(-\alpha\log(x))$ (as $p_j(\cdot)$ is bounded, continuous),
the third term vanishes $J_3(n) \rightarrow  0$. \qed

\subsection{Power-laws for limiting empirical degree distribution}
Recall, with respect to the definition of $\pi_j$ (\ref{pi_j}),
that $\{D(y): y\geq 0\}$
is a Markov branching process with exponential($\lambda=1$) lifetime
distribution, offspring distribution $\{p'_j=p_{j-1}\}_{j\ge 2}$,
immigration rate $\beta\geq 0$, immigration size distribution
$\{p_j\}_{j\geq 1}$, and initial value $D(0)$ distributed according
to $\{p_j\}_{j\ge 1}$.

\vskip .2cm

\noindent {\bf Proof of Theorem \ref{powerlaw}.}  First note, for
$s\ge 0$, that
\begin{eqnarray}
\sum_{j\geq1}j^s \pi_j & = & (2m+\beta)\int_0^\infty
e^{-(2m+\beta)y} \sum_{j\geq 1}j^sp_j(y) dy \nonumber\\
& =&
 (2m+\beta)\int_0^\infty
e^{-(2m+\beta)y} E[(D(y))^s] dy.\label{pi_j_mean}
\end{eqnarray}
By Proposition \ref{prop5}, $D(y)e^{-my} \rightarrow \zeta$ a.s., and
when $\sum_{j\geq 1}j^{2+\beta/m}p_j<\infty$, we have by
Proposition \ref{appendix1} that
$E[(\sup_{y\geq 0} D(y)e^{-my})^{2+\beta/m}]<\infty$.  Hence, by
dominated convergence, for $s\leq 2 + \beta/m$, and a constant $C$,
$$ Ce^{msy} \ \leq E[(D(y))^s] \ \leq C^{-1}e^{msy}.$$
Plugging into (\ref{pi_j_mean}), we prove the theorem for all
$0\leq s\leq 2+\beta/m$.  Noting $\sum_{j\geq 1}j^{s}\pi_j\geq \sum_{j\geq 1}
j^{2+\beta/m}\pi_j$ for $s\geq 2+\beta/m$ finishes the proof. \qed

We now evaluate $\pi_j$ in a special  case.  The formula, similar to
that in \cite{Mori1} and section 4.2 \cite{RTV}, gives the
asymptotics mentioned in Remark \ref{power_law_rm}.

\begin{prop}
\label{explicit_pi}
When $X_i\equiv x_0$ for an integer $x_0\geq 1$, we have for $j\geq 1$
that
$$\pi_j  \ = \ \frac{2x_0+\beta}{j+2x_0 +2\beta}\prod_{k=1}^{l-1}
\frac{kx_0+\beta}{(k+2)x_0 +2\beta}$$ when $j = lx_0$ for $l\geq 1$
(where the product is set equal to $1$ when $l=1$), and $\pi_j = 0$
otherwise. Hence, for large $j$,  $\pi_j = O(j^{-[3+\beta/x_0]})$ when $j = lx_0$
for $l\geq 1$, and $\pi_j = 0$ otherwise.
\end{prop}

{\it Proof.}  First note as $X_i\equiv x_0$ that $m=D(0)=x_0$, and the
process $D(\cdot)$ moves in steps of $x_0$.  Clearly, then $\pi_j=0$
when $j$ is not a multiple of $x_0$.  When $j=lx_0$ for $l\geq 1$, let $A_j$ be the first time the
process
$D(\cdot)$ equals $j$,
$$A_j \ = \ \inf \big\{ y\geq 0: D(y) = j\big\} \ < \ \infty \ \ \
{\rm a.s.}$$
Let $B_j$ be the time, after $A_j$, that the process
$D(\cdot)$ spends at $j$; note that conditioned on $A_j$, $B_j$ is an
exponential$(j+\beta)$ variable.  Then, we write
\begin{eqnarray*}
\pi_j &=& (2x_0+\beta)\int_0^\infty e^{-y(2x_0+\beta)} p_j(y) dy\\
&=& (2x_0+\beta) E\bigg[\int_{A_j}^{A_j +B_j}
e^{-y(2x_0+\beta)} dy \bigg]\\
&=&
E\bigg[e^{-A_j(2x_0+\beta)}\big[1-e^{-B_j(2x_0+\beta)}\big]\bigg]\\
&=& E[e^{-A_j(2x_0+\beta)}]\bigg[1-\frac{j+\beta}{j+2x_0+
    2\beta}\bigg]\\
&=&\frac{2x_0+\beta}{j+2x_0
  +2\beta}E[e^{-A_j(2x_0+\beta)}].
\end{eqnarray*}

As $X_i\equiv x_0$, for $j= lx_0$ and $l\geq 2$, we have
$A_j$ is the sum of $l-1$ independent exponential random variables with
parameters
$x_0+\beta,\ldots,(l-1)x_0 +\beta$, and so
$$E[e^{-A_j(2x_0+\beta)}] \  = \ \prod_{k=1}^{l-1}
\frac{kx_0+\beta}{(k+2)x_0 +2\beta}.$$
When $l=1$, then $A_j =0$, giving $\pi_{x_0} = (2x_0 +\beta)/(3x_0 +2\beta)$. \qed

\bibliographystyle{plain}

\begin{thebibliography}{99}
\frenchspacing

\addcontentsline{toc}{chapter}
{\protect\numberline{12.}{Bibliography}}

\bibitem{Albert_Barabasi}
Albert, R., Barabasi, A.-L. (2002) Statistical mechanics of
complex networks. {\it Reviews of Modern Phys.} {\bf 74} 47-97.


\bibitem{Asmussen}
Asmussen, S., Hering, H. (1983) {\it Branching Processes.}
 Progress in Probability and Statistics {\bf 3}, Birkh\"auser, Boston.

\bibitem{athreya}
Athreya, K.B. (1969) On the equivalence of conditions on a
branching process in continuous time and on its offspring
distribution. {\it J. Math. Kyoto Univ.} {\bf 9} 41-53.



\bibitem{athreya-karlin}
Athreya, K.B., Karlin S. (1967) Limit theorems for the split times
of branching processes. {\it J. Math. Mech.} {\bf 17} 257-277.




\bibitem{AN}
Athreya, K.B., Ney, P. (2004) {\it Branching Processes.} Dover, New York.

\bibitem{AL}
Athreya, K.B., Lahiri, S.N.  (2006) {\it Measure Theory and
Probability Theory.} Springer, New York.


\bibitem{BA}
Barabasi, A.-L., Albert, R. (1999) Emergence of scaling in random networks. {\it
  Science} {\bf 286} 509-512.

\bibitem{Bollobas}
Bollobas, B., Riordan, O, Spencer, J., Tusnady, G. (2001) The degree
sequence of a scale-free random graph process. {\it Random Structures
  and Algorithms} {\bf 18} 279-290.


\bibitem{CL}
Chung, F., Lu, L. (2006) {\it Complex Graphs and Networks} CBMS {\bf
  107}, American Mathematical Society,
Providence.

\bibitem{CF}
Cooper, C., Frieze, A. (2003) A general model of web graphs.
{\it Random Structures Algorithms} {\bf 22} 311-335.


\bibitem{Durrett_book}
Durrett, R. (2006) {\it Random Graph Dynamics} Cambridge University
Press, Boston.

\bibitem{Mitzenmacher}
Mitzenmacher, M. (2006) A brief history of generative models for power law and
lognormal distributions. {\it Internet Math.} {\bf 1} 226-251.


 \bibitem{Mori1} Mori, T. (2002) On random trees. {\it Studia
 Sci. Math. Hungar.} {\bf 39} 143-155.

\bibitem{Mori}
Mori, T. (2005) The maximum degree of the Barabasi-Albert random
tree. {\it Comb. Probab. Computing}, {\bf 14}, 339-348.

\bibitem{Newman} Newman, M.E.J. (2003) The structure and function of
  complex networks. {\it SIAM Review} {\bf 45} 167-256.

 \bibitem{OS}
 Oliveira, R., Spencer, J. (2006) Connectivity transitions in networks
 with super-linear preferential attachment. {\it Internet Math.} {\bf
 2} 121-163.

\bibitem{Pemantle}
Pemantle, R. (2005) Random processes with reinforcement.  Pre-print
available at http://www.math.upenn.edu/~pemantle/papers/Papers.html

\bibitem{RTV}
Rudas, A., Toth, B., Valko, B. (2006) Random trees and general
branching processes. {\it to appear Random Structures and Algorithms,}
{\it arXiv.math.PR/0503728 v2}.

\bibitem{Simon} Simon, H.A. (1955) On a class of skew distribution
  functions. {\it Biometrika} {\bf 42} 425-440.

\bibitem{Yule} Yule, G. (1925) A mathematical theory of evolution
  based on the conclusions of Dr. J.C. Willis.  {\it
  Phil. Trans. Roy. Soc. Lond. Ser. B} {\bf 213} 21-87.



\end{thebibliography}

\end{document}